\documentclass [11pt]{article}
\usepackage [a4paper]{geometry}
\usepackage [T1]{fontenc}
\usepackage {ae}
\usepackage [utf8]{inputenc}
\usepackage {amsmath,amssymb,amsfonts}
\usepackage {hyperref}

\newtheorem {theorem}{Theorem}
\newtheorem {corollary}[theorem]{Corollary}
\newenvironment {proofof}[1]
   {{\medskip \parindent0cm \bf Proof#1:}}
   {\hspace* {\fill} $\Box$}
\newenvironment {proof}
   {\begin {proofof}{}}
   {\end {proofof}}

\newcommand {\af}{\mathfrak a}
\newcommand {\cf}{\mathfrak c}
\newcommand {\kf}{\mathfrak k}
\newcommand {\nf}{\mathfrak n}
\newcommand {\pf}{\mathfrak p}
\newcommand {\zf}{\mathfrak z}
\newcommand {\w}[1]{\mathfrak {w}_{#1}}
\renewcommand {\wp}{\w {p_1}}
\newcommand {\wq}{\w {p_2}}
\newcommand {\wN}{\w {N}}
\newcommand {\wpq}{\w {p_1,p_2}}
\newcommand {\wpk}{\w {p_1, \ldots, p_k}}
\newcommand {\alN}{|_N}
\newcommand {\Z}{\mathbb Z}
\newcommand {\Q}{\mathbb Q}
\newcommand {\R}{\mathbb R}
\newcommand {\C}{\mathbb C}
\newcommand {\F}{\mathbb F}
\newcommand {\Fc}{\mathcal F}

\newcommand {\Fn}{\Fc_n}
\newcommand {\Oc}{\mathcal O}
\newcommand {\OD}{\Oc_D}
\newcommand {\legendre}[2]{\genfrac {(}{)}{1pt}{}{#1}{#2}}
\newcommand {\Omegac}{\Omega_c}
\newcommand {\Cl}{\operatorname {Cl}}
\newcommand {\Sl}{\operatorname {Sl}}
\newcommand {\Gal}{\operatorname {Gal}}

\renewcommand {\epsilon}{\varepsilon}
\renewcommand {\geq}{\geqslant}
\renewcommand {\leq}{\leqslant}

\title {Singular values of multiple eta-quotients for ramified primes}
\author {Andreas Enge\footnote {INRIA, LFANT, F-33400 Talence, France \newline
CNRS, IMB, UMR 5251, F-33400 Talence, France \newline
Univ. Bordeaux, IMB, UMR 5251, F-33400 Talence, France \newline
andreas.enge@inria.fr \newline
This research was partially funded by ERC Starting Grant ANTICS 278537.}
\ and
Reinhard Schertz\footnote {Universität Augsburg \newline
schertz@math.uni-augsburg.de}}
\date {6 May 2013}

\begin {document}
\maketitle

\begin {abstract}
We determine the conditions under which singular values of multiple
$\eta$-quotients of square-free level, not necessarily prime to~$6$, yield
class invariants, that is, algebraic numbers in ring class fields of
imaginary-quadratic number fields. We show that the singular values lie in
subfields of the ring class fields of index $2^{k' - 1}$ when $k' \geq 2$
primes dividing the level are ramified in the imaginary-quadratic field,
which leads to faster computations of elliptic curves with prescribed
complex multiplication. The result is generalised to singular values of
modular functions on $X_0^+ (p)$ for $p$ prime and ramified.
\end {abstract}


\bigskip

Let $K= \Q (\sqrt \Delta)$ be an imaginary-quadratic number field of
discriminant $\Delta$, and let $\OD$ be the order of discriminant $D = c^2
\Delta$ and conductor $c$ in $K$. For a modular function $f$ and an argument
$\tau \in K \subseteq \C$ with $\Im \tau > 0$ we call the \textit {singular
value} $f (\tau)$ a \textit {class invariant} if it lies in the \textit
{ring class field} $\Omegac$.
This is the abelian extension of $K$ with Galois group
canonically isomorphic to $\Cl (\OD)$ through the
\textit {Artin map} $\sigma : \Cl (\OD) \to \Gal (\Omegac / K)$,
which sends a prime ideal representing an ideal class to its associated
\textit {Frobenius automorphism}.
In this article, we are interested in class
invariants derived from multiple $\eta$-quotients, and we examine in
particular cases where those generate a subfield of the ring class field.

In \S1 we define the multiple $\eta$-quotients under consideration and
collect their properties, in particular their transformation behaviour under
unimodular matrices. We then proceed in \S2 to determine conditions under
which their singular values lie in the ring class field and show how to
compute their characteristic polynomials with respect to the number field
extension $\Omegac / K$ using roots of an \textit {$n$-system}, suitably
normalised quadratic forms of discriminant~$D$ representing the class group.
The results of \S3 are at the heart of this article: We show that if some or
all of the primes dividing the level of the multiple $\eta$-quotient are
ramified in~$K$, then the singular values lie in fact in a subfield $L$ of
$\Omega_c$ of index a power of~$2$; more precisely, the Galois group of
$\Omegac / L$ is elementary abelian, so that $\Omegac$ is a compositum of
linearly disjoint quadratic extensions of~$L$. An alternative proof for the
special case that all primes are ramified leads to an 
interesting generalisation in \S4, namely to functions of prime level
invariant under the Fricke--Atkin--Lehner involution. We conclude by some
class field theoretic remarks in \S5, showing how this computationally less
expensive construction of subfields of $\Omegac$ may be completed to obtain
the full ring class field.

\section {Multiple $\eta$-quotients}
\label {sec:functions}

In this section, we define multiple $\eta$-quotients and collect their basic
properties, most of which are either well-known or readily verified.

Let $\eta$ be Dedekind's function, and consider positive integers $p_1,
\ldots, p_k$; in later sections, they will be distinct primes, but this
restriction is not needed for the time being.

The \textit {simple $\eta$-quotient} of level $p_1$ is defined by
\begin {equation}
\label {eq:simpleeta}
\wp (z) = \frac {\eta \left( \frac {z}{p_1} \right)}{\eta (z)},
\end {equation}
and the \textit {double $\eta$-quotient} of level $p_1 p_2$ by
\begin {equation}
\label {eq:doubleeta}
\wpq (z) = \frac {\eta \left(\frac {z}{p_1} \right)
\eta \left(\frac {z}{p_2} \right)}
{\eta \left( \frac {z}{p_1 p_2} \right) \eta (z)}
=
\frac {\wp (z)}{\wp \left( \frac {z}{p_2} \right)}
=
\frac {\wq (z)}{\wq \left( \frac {z}{p_1} \right)}
.
\end {equation}
The process may be continued inductively by letting
\begin {equation}
\label {eq:multieta}
\w {p_1, \ldots, p_{k+1}}
= \frac {\w {p_1, \ldots, p_k} (z)}
{\w {p_1, \ldots, p_k} \left( \frac {z}{p_{k+1}} \right)},
\end {equation}
so that $\wpk$ is a quotient of transformed $\eta$-functions with $2^k$
factors in the numerator and as many in the denominator.

Let $\Fn$ denote the set of modular functions of level~$n$ whose
$q$-expansions have coefficients in $\Q (\zeta_n)$. The powers of the
multiple $\eta$-quotients are elements of $\Fn$ for some~$n$, which can be
determined from the transformation behaviour of $\eta$ under unimodular
substitutions. Here and in the following, we consider unimodular matrices $M
= \begin {pmatrix} a & b \\ c & d \end {pmatrix} \in \Gamma = \Sl_2 (\Z) /
\{ \pm 1 \}$, normalised such that $c \geq 0$, and $d > 0$ if $c = 0$. We
denote by $\cdot'$ the odd part and by $\lambda (\cdot)$ the $2$-adic
valuation of a number, so that $c = c' \, 2^{\lambda (c)}$; by convention,
$\lambda (0) = 0' = 1$. Let
\begin {eqnarray*}
\label {eq:ebar}
\bar e (M) & = &  a b + c \big( d (1 - a^2) - a \big) + 3 c' (a-1) \\
\label {eq:e}
e (M) & = & \bar e (M) + \frac {3}{2} \lambda (c) (a^2 -1) \\
\label {eq:epsilon}
\epsilon (M) & = & \legendre {a}{c'} \zeta_{24}^{e (M)}
\end {eqnarray*}
with $\zeta_{24} = e^{2 \pi i / 24}$. Notice that $e (M) - \bar e (M)$ is
divisible by~$12$, since $a$ and $c$ cannot be even simultaneously. By
\cite[Th.~1.10.1]{Schertz10},
\[
\eta (M z) = \epsilon (M) \sqrt {c z + d} \, \eta (z).
\]
So for $M \in \Gamma^0 (N)$ with $b = N b_0$ we have, cf.
\cite[Th.~3]{EnSc05},
\begin {eqnarray}
\label {eq:trafosimpleeta}
\wN (M z) & = & \epsilon \begin {pmatrix} a & b_0 \\ N c & d \end {pmatrix}
\epsilon \begin {pmatrix} a & N b_0 \\ c & d \end {pmatrix}^{-1} \wN (z) \\
\nonumber
& = & \legendre {a}{N'} \zeta_{24}^{(N-1) \big( -a b_0 + c (d (1 - a^2) -a)
\big) + 3 (N' - 1) c' (a-1) + \frac {3}{2} (\lambda (N c) - \lambda (c))
(a^2 - 1)} \wN (z),
\end {eqnarray}
where $\lambda (N c) - \lambda (c)$ equals $0$ for $c = 0$, and
$\lambda (N)$ otherwise.
With $s = \frac {24}{\gcd (24, N-1)}$ and $e \mid s$, we have $\wN^e \in
\Fc_{\frac {s}{e} N}$ by \cite[Th.~6]{EnMo09}.

For the double $\eta$-quotients $\wpq$ of level $N = p_1 p_2$ and $M \in
\Gamma^0 (N)$ with $b = N b_0$ we compute
\begin {eqnarray}
\label {eq:trafodoubleeta}
\wpq (Mz) & = & \epsilon
\begin {pmatrix} a & p_2 b_0 \\ p_1 c & d \end {pmatrix}
\epsilon \begin {pmatrix} a & p_1 b_0 \\ p_2 c & d \end {pmatrix}
\epsilon \begin {pmatrix} a & b_0 \\ N c & d \end {pmatrix}^{-1} \epsilon
\begin {pmatrix} a & N b_0 \\ c & d \end {pmatrix}^{-1} \wpq (z) \\
\nonumber
& = & \zeta_{24}^{-(p_1 - 1)(p_2 - 1) \big( a b_0 + c (d (1 - a^2) -a) \big)
- 3 (p_1' - 1)(p_2' - 1) c' (a-1)} \wpq (z).
\end {eqnarray}
Let $s = \frac {24}{\gcd (24, (p_1 - 1)(p_2 - 1))}$ and $e \mid s$. Then
$\wpq^e \in \Fc_{\frac {s}{e} N}$, cf.~\cite [Th.~7]{EnSc05} and \eqref
{eq:trafodoubleeta}.

For $\eta$-quotients of order~$k \geq 2$, where $N = p_1 \cdots p_k$, the
formula generalises as
\begin {equation}
\label {eq:trafomultieta}
\begin {split}
& \wpk (M z) = \\
& \zeta_{24}^{-(p_1 - 1) \cdots (p_k - 1) \left( a b_0 + (-1)^k c (d (1 -
a^2) - a) \right) - 3 (-1)^k (p_1' - 1) \cdots (p_k' - 1) c' (a - 1)} \wpk
(z).
\end {split}
\end {equation}
If $s = \frac {24}{\gcd (24, (p_1 - 1) \cdots (p_k - 1))}$ and $e \mid s$,
then $\wpk^e \in \Fc_{\frac {s}{e} N}.$

For later reference, we also recall the transformation behaviour of
$\gamma_2 = \sqrt [3]{j}$ and $\gamma_3 = \sqrt {j - 1728}$, see
\cite[\S2.4.3]{Schertz10}:
\begin {equation}
\label {eq:trafogamma}
\gamma_2 (M z) = \zeta_3^{- e (M)} \gamma_2 (z),
\gamma_3 (M z) = (-1)^{e (M)} \gamma_3 (z),
(\gamma_2 \gamma_3) (M z) = \zeta_6^{e (M)} (\gamma_2 \gamma_3) (z).
\end {equation}

Since $\eta$ has a rational $q$-expansion, so does $\wpk$. For $S = \begin
{pmatrix} 0 & -1 \\ 1 & 0 \end {pmatrix}$, we have $\eta \left( \frac
{Sz}{N} \right) = \epsilon (S) \sqrt {N z} \, \eta (N z)$ by
\cite[Th.~3]{EnSc05}, so that
\[
\wN (Sz) = \sqrt N \, \frac {\eta (N z)}{\eta (z)}
= \sqrt N \, \wN (N z)^{-1}
\]
has a rational $q$-expansion up to a
factor $\sqrt N$. For $k \geq 2$ and $N = p_1 \cdots p_k$, this implies
$\wpk (S z) = \wpk (N z)^{(-1)^k}$, which has a rational $q$-expansion.

Denote by $\alN$ the Fricke-Atkin-Lehner involution associated to $\Gamma^0
(N)$, so that $f\alN (z) = f \left( \frac {-N}{z} \right)$. The previous
equation can be rewritten as
\begin {equation}
\label {eq:fal}
\wpk\alN (z) = \wpk (z)^{(-1)^k};
\end {equation}
in particular, $\wpk$ is invariant under the involution for even $k$.

\section {Singular values of multiple $\eta$-quotients}

\subsection {Class invariants}
\label {ssec:classinv}

A very general result on class invariants is obtained in
\cite[Th.~4]{Schertz02}: Let $f$ be modular for $\Gamma^0 (N)$ such that $f$
and $f \circ S$ have rational $q$-expansions. Assume that there are $A$ and
$B$ such that $\af = A \Z +  \frac {-B + \sqrt D}{2} \Z$ is a proper ideal
of $\OD$, $N \mid \frac {B^2 - D}{4 A}$ and $N$ is coprime to the
conductor~$c$. Let $\tau = \frac {-B + \sqrt D}{2 A}$ be the basis quotient
of $\af$ with $\Im \tau > 0$. Then $f (\tau) \in \Omegac$. We use this
theorem to determine when singular values of powers of multiple
$\eta$-quotients are class invariants.

\begin {theorem}
\label {th:singulardouble}
Let $p_1 < p_2$ be primes and $e$ an integer such that one of the following
conditions is satisfied:
\begin {enumerate}
\item
$\{ p_1, p_2 \} \cap  \{ 2, 3 \} = \emptyset$ \\
$e (p_1 - 1)(p_2 - 1) \equiv 0 \pmod 3$ or  $3 \nmid D$, and \\
$e (p_1 - 1)(p_2 - 1) \equiv 0 \pmod 8$ or  $2 \nmid D$.
\item
$p_1 = 3$, $p_2 \geq 5$ \\
$e (p_2 - 1) \equiv 0 \pmod 3$, and \\
$e (p_2 - 1) \equiv 0 \pmod 4$ or $2 \nmid D$.
\item
$p_1 = 2$ \\
$e (p_2 - 1) \equiv 0 \pmod {24}$.
\end {enumerate}
Let $p_1$, $p_2$ be non-inert in $K$ and not dividing the conductor~$c$, and
let $\tau = \frac {-B + \sqrt D}{2 A}$ be such that $p_1 p_2 \mid C = \frac
{B^2 - D}{4 A}$.
Then $\wpq^e (\tau) \in \Omegac$.
\end {theorem}

\begin {proof}
The result is well-known in the case $e (p_1 - 1)(p_2 - 1) \equiv 0 \pmod
{24}$, cf.~\cite[Th.~3.2]{EnSc04}, since then $\wpq^e$ is modular for
$\Gamma^0 (p_1 p_2)$.
In all other cases, $p_1$ and $p_2$ are odd, and \eqref {eq:trafodoubleeta}
becomes
\begin {equation}
\label {eq:trafodoubleetaodd}
\wpq (M z) = \zeta_6^{- \frac {(p_1 - 1)(p_2 - 1)}{4} \left( e (M) + a b_0
(1 - p_1 p_2) \right)} \wpq (z).
\end {equation}
Consider first the case that $p_1, p_2 \neq 3$, and let $h = (\gamma_2
\gamma_3)^{\frac {(p_1 - 1)(p_2 - 1)}{4}} \wpq$. From \eqref
{eq:trafogamma}, \eqref {eq:trafodoubleetaodd} and the fact that at least
one of $p_1 - 1$, $p_2 - 1$ and $1 - p_1 p_2$ is divisible by~$3$, we deduce
that $h$ is modular for $\Gamma^0 (p_1 p_2)$, cf.~\cite[\S2.4.3]{Schertz10},
so that $h (\tau) \in \Omegac$. Since $\gamma_2 (\tau) \in \Omegac$ for $3
\nmid D$ and $\gamma_3 (\tau) \in \Omegac$ for $2 \nmid D$ by \cite[Th.~2
and~3]{Schertz02}, and $j (\tau) = \gamma_2 (\tau)^3 = \gamma_3 (\tau)^2 +
1728 \in \Omegac$, we have $(\gamma_2 \gamma_3)^{e \frac {(p_1 - 1)(p_2 -
1)}{4}} (\tau) \in \Omegac$ under the assumptions of the theorem, which
proves the result.

If $p_1 = 3$, the function $h^e$ is modular for $\Gamma^0 (p_1 p_2)$ under
the conditions of the theorem, and we conclude analogously.
\end {proof}

Similar results hold for higher order $\eta$-quotients; in fact, adding odd
primes makes it easier to satisfy the restrictions modulo~$8$.

\begin {theorem}
\label {th:singularmulti}
For $k \geq 3$, let $e$ be an integer and $p_1, \ldots, p_k$ distinct primes
such that one of the following conditions is satisfied:
\begin {enumerate}
\item
$e (p_1 - 1) \cdots (p_k - 1) \equiv 0 \pmod {24}$.
\item
All $p_i$ are odd and congruent to~$-1$ modulo~$3$, and $3 \nmid D$.
\end {enumerate}
Let $p_1, \ldots, p_k$ be non-inert in $K$ and not dividing the
conductor~$c$, and let $\tau = \frac {-B + \sqrt D}{2 A}$ be such that $N =
p_1 \cdots p_k \mid C = \frac {B^2 - D}{4 A}$.
Then $\wpk^e (\tau) \in \Omegac$.
\end {theorem}

\begin {proof}
The first case is trivial. In the second case, we use the auxiliary function
$h = \gamma_2^{(-1)^k} \wpk$, which is modular for $\Gamma^0 (N)$ by \eqref
{eq:trafomultieta}, \eqref {eq:trafogamma}, where $N \equiv (-1)^k \pmod 3$,
and we conclude as in the proof of Theorem~\ref {th:singulardouble}.
\end {proof}

\subsection {$n$-systems and reality of class polynomials}

In order to compute the \textit {class polynomial} $H_D^f (X)$, the
characteristic polynomial of a class invariant $f (\tau)$ under the action
of $\Cl (\OD)$, we need to explicitly determine the conjugates of $f (\tau)$
under the Galois group $\Cl (\OD)$. Classically, this is done using Shimura
reciprocity. The concept of $n$-systems was introduced in \cite{Schertz02};
it encapsulates Shimura reciprocity and provides a simple way of expressing
the conjugates as singular values of the \textit {same} function $f$ in
basis quotients of suitably normalised quadratic ideals.

An $n$-system for $\Cl (\OD)$ is defined as a set of quadratic forms $[A_i,
B_i, C_i] = A_i X^2 + B_i X + C_i$, $C_i = \frac {B_i^2 - D}{4 A_i}$, $\gcd
(A_i, B_i, C_i) = 1$ for $1 \leq i \leq h (\OD)$, such that the ideals
$\af_i = \begin {pmatrix} \frac {-B_i + \sqrt D}{2} \\ A_i \end
{pmatrix}_\Z$ form a system of representatives of $\Cl (\OD)$, and
furthermore
\[
\gcd (A_i, n) = 1 \text { and } B_i \equiv B_1 \pmod {2 n}.
\]
(Here and in the following, we write $\Z$-bases of ideals as column vectors
such that the \textit {basis quotient} of its two entries is the root of the
quadratic form with positive imaginary part.)

Notice that if $n \mid C_1$, then $n \mid C_i$ for all~$i$. For any~$n$, an
initial $[A_1, B_1, C_1]$ with $\gcd (A_1, n) = 1$ is easily extended to an
$n$-system using the algorithm of \cite[Prop.~3]{Schertz02}. It is shown in
\cite[Th.~7]{Schertz02} that if $f \in \Fn$ is such that $f \circ S$ has a
rational $q$-expansion, $\tau_i = \frac {-B_i + \sqrt D}{2 A_i}$ are basis
quotients coming from an $n$-system and $f (\tau_1) \in \Omegac$, then its
algebraic conjugates are the $f (\tau_i)$. More precisely, if
$\sigma : \Cl (\OD) \stackrel{\simeq}{\to} \Gal (\Omegac / K)$
is the Artin map, then
\[
H_D^f (X) = \prod_{\kf \in \Cl (\OD)} \big( X - f (\tau_1)^{\sigma (\kf)}
\big)
= \prod_{i=1}^{h (\OD)} \big( X - f (\tau_i) \big)
\in K [X].
\]

This characterisation can also be used to identify pairs of complex
conjugate roots of the class polynomial whenever the latter is real. We
recall that $\alN$ denotes the Fricke--Atkin--Lehner involution such that
$f\alN (z) = f \left( \frac {-N}{z} \right)$.

\begin {theorem}
\label {th:realnsystem}
Let $N$ and $n$ be integers such that $N \mid n$, let $f \in \Fn$ be such
that $f$ and $f \circ S$ have rational $q$-expansions, and assume that there
is an ideal $\nf = \begin {pmatrix} \frac {-B_1 + \sqrt D}{2} \\ A_1 \end
{pmatrix}$ of basis quotient $\tau_1 = \frac {-B_1 + \sqrt D}{2 A_1}$ and
cofactor $C_1 = \frac {B_1^2 - D}{4 A_1}$ such that $C_1 = N$ and $f
(\tau_1) \in \Omegac$.

Write $\tau_0 = A_1 \tau_1$. Then $f (\tau_0) = f (\tau_1)^{\sigma (\nf)}
\in \Omegac$.

Assume that there is a proper ideal $\cf$ of~$\OD$ such that $f\alN (\tau_0)
= f (\tau_0)^{\sigma (\cf)}$. Then $H_D^f (X) \in \Q [X]$. More precisely,
if the $\af_i$ are given by an $n$-system $[A_i, B_i, C_i]$ for $\Cl (\OD)$
with basis quotients $\tau_i$ and $\sim$ denotes equivalence in the class
group, then $f (\tau_i)$ and $f (\tau_j)$ are complex conjugates if $\af_i
\af_j \sim \nf \cf^{-1}$. In particular, $f (\tau_i) \in \R$ if $\af_i^2
\sim \nf \cf^{-1}$.
\end {theorem}

This result generalises \cite[Th.~3.4]{EnSc04}, which treats the case that
$f\alN = f$ (so that $\cf = \OD$) and that $N = n$. The latter condition is
used in \cite{EnSc04} only to ensure that $f (\tau_1) \in \Omegac$, which
instead we added to the hypotheses of the theorem. Once the correct
assumptions are identified, the proof itself is very similar.

\begin {proof}
The proof of \cite[Th.~7]{Schertz02} shows that
\begin {equation}
\label {eq:tau0}
f (\tau_i) = f (\tau_0)^{\sigma (\af_i^{-1})}
\end {equation}
for all $i$. (Here we need the rationality of the $q$-expansion of $f \circ
S$, which implies that all $g_i$ equal $g = f$ in the notation of
\cite{Schertz02}.)

Denote by $\kappa$ the complex conjugation, and recall that $\Gal (\Omegac /
\Q)$ is isomorphic to the generalised dihedral group $\Cl (\OD) \rtimes
\langle \kappa \rangle$ with $\kappa \sigma (\af) \kappa = \sigma
(\af^{-1})$ for $\af \in \Cl (\OD)$.

We first consider $i = 1$ and use that $f (z)^\kappa = f (- z^\kappa)$ by
the rationality of the $q$-expansion of $f$. Then
\begin {eqnarray*}
f (\tau_1)^\kappa
& = & f \left( \frac {-B_1 + \sqrt D}{2 A_1} \right)^\kappa
= f \left( \frac {B_1 + \sqrt D}{2 A_1} \right)
= f \left( \frac {2 N}{B_1 - \sqrt D} \right)
\text { since $C_1 = N$} \\
& = & f\alN (\tau_0)
= f (\tau_0)^{\sigma (\cf)}
\text { by assumption.}
\end {eqnarray*}
Let now $i$ and $j$ be such that $\af_i \af_j \sim \nf \cf^{-1}$. We compute
\begin {eqnarray*}
f (\tau_i)^\kappa
& = & f (\tau_1)^{\sigma (\nf \af_i^{-1}) \kappa}
\text { by \eqref {eq:tau0}} \\
& = & f (\tau_1)^{\kappa \sigma (\nf^{-1} \af_i)}
= f (\tau_0)^{\sigma (\cf \nf^{-1} \af_i)}
= f (\tau_j)^{\sigma (\af_j \cf \nf^{-1} \af_i)}
= f (\tau_j).
\end {eqnarray*}
\end {proof}

\begin {corollary}
\label {cor:realeven}
\begin {sloppypar}
Under the hypotheses of Theorems~\ref {th:singulardouble} or~\ref
{th:singularmulti}, let $N = p_1 \cdots p_k$, $s = \frac {24}{\gcd \big( 24,
(p_1 - 1) \cdots (p_k - 1)) \big)}$, $e \mid s$ and $n = \frac {s}{e} N$.
Then there is an $n$-system satisfying $C_1 = N$, which yields the roots of
the class polynomial $H_D^{\wpk^e}$.
\end {sloppypar}

Moreover, the polynomial is real for even $k$; its complex conjugate roots
may be identified by Theorem~\ref {th:realnsystem} with $\cf = \OD$.
\end {corollary}

\begin {proof}
Since none of the prime divisors of~$N$ is inert, there is an ideal of
norm~$N$ or, equivalently, a quadratic form $[A_1, B_1, C_1]$ with $C_1 =
N$. An application of Theorem~\ref {th:realnsystem} using \eqref {eq:fal}
finishes the proof.
\end {proof}

For $f = \wpk^e$ with odd $k$, the class polynomial is in general defined
over $K$. As for simple $\eta$-quotients in \cite[Th.~21]{EnMo09}, a
particular case is easily identified in which the class polynomial is real:
If $n \mid B_1$, which implies that $N \mid D$ and that all primes dividing
$N$ are ramified, then $- B_i \equiv B_i \pmod {2 n}$, and for every ideal
class $\af_i$ represented by an element $[A_i, B_i, C_i]$ of the $n$-system,
the inverse class $\overline \af_i$ is represented by the element $[A_i,
-B_i, C_i]$. Together with the rationality of the $q$-expansion of $f$, this
implies that $H_D^f$ is real and that $f (\tau_j) = \overline {f (\tau_i)}$
if $\af_j \sim \overline \af_i \sim \af_i^{-1}$. It is then enough to
compute only $\frac {h (\OD) + h_0}{2}$ values $f (\tau_i)$, where $h_0$ is
the number of real roots of the class polynomial, which is bounded above by
the size of the $2$-torsion subgroup of $\Cl (\OD)$. So it is generally
small, and the required number of function evaluations
boils down to essentially $\frac {h (\OD)}{2}$.

We show in Corollary~\ref {cor:realodd} that the condition $n \mid B_1$ is
in fact too restrictive: $H_D^f$ is real already when only one prime
dividing $N$ is ramified in $K$.

\subsection {Examples}

Let $D = -215 = -5 \cdot 31$, for which $2$, $3$, $7$ and $11$ split and $5$
ramifies. The class number of $\Oc_{-215}$ is $14$, and $\Oc_{-215} = \Z +
\omega \Z$ with $\omega = \frac {1 + \sqrt {-215}}{2}$.

Using the double $\eta$-quotient for the primes $7$ and $11$, the full
exponent $s$ equals~$2$, but we may use the lower exponent $e=1$ by
Theorem~\ref {th:singulardouble}:
\[
\begin {split}
H_{-215}^{\w {7, 11}} (X) = &
X^{14} -10 X^{13} + 42 X^{12} -97 X^{11} + 144 X^{10} -147 X^9 + 89 X^8 + 25
X^7 \\
& -124 X^6 + 113 X^5 -23 X^4 -28 X^3 + 20 X^2 -5 X + 1.
\end {split}
\]
The triple $\eta$-quotient for $2$, $3$ and $7$ has $s = e = 2$, and
\[
\begin {split}
H_{-215}^{\w {2, 3, 7}^2} (X) = &
X^{14} + (17+\omega) X^{13} + (104+16 \omega) X^{12} + (211+107 \omega)
X^{11} \\
& + (-573+379 \omega) X^{10} + (-4197+737 \omega) X^9 + (-10230+686 \omega)
X^8 \\
& -13247 X^7 + (-9544-686 \omega) X^6 + (-3460-737 \omega) X^5 \\
& + (-194-379 \omega) X^4 + (318-107 \omega) X^3 + (120-16 \omega) X^2 \\
& + (18-\omega) X + 1.
\end {split}
\]
With the ramified $5$ instead of the split $7$, the polynomial becomes real,
but one needs the higher exponent $s = e = 3$:
\[
\begin {split}
H_{-215}^{\w {2, 3, 5}^3} (X) = &
X^{14} + 22 X^{13} + 175 X^{12} + 578 X^{11} + 819 X^{10} + 2190 X^9 + 10130
X^8 \\
& + 17295 X^7  + 10130 X^6 + 2190 X^5 + 819 X^4 + 578 X^3 + 175 X^2 + 22 X +
1.
\end {split}
\]
Examples for $\eta$-quotients of order~$4$ and~$5$ are given by
\[
\begin {split}
H_{-215}^{\w {2, 3, 5, 7}} (X) = &
X^{14} - X^{13} -8 X^{12} -12 X^{11} -7 X^{10} -4 X^9 -17 X^8 -29 X^7 -17
X^6 \\
& -4 X^5 -7 X^4 -12 X^3 -8 X^2 - X + 1
\end {split}
\]
and
\[
\begin {split}
H_{-215}^{\w {2, 3, 5, 7, 11}} (X) = &
X^{14} -3 X^{13} + 6 X^{12} + 35 X^{11} + 80 X^{10} + 130 X^9 + 188 X^8 +
201 X^7 \\
& + 188 X^6 + 130 X^5 + 80 X^4 + 35 X^3 + 6 X^2 -3 X + 1.
\end {split}
\]

\section {Singular values for ramified primes}
\label {sec:ramified}

\subsection {Class invariants in subfields}

In this section, we show that the singular values of multiple
$\eta$-quotients lie in subfields of the ring class field when at least two
of the involved primes are ramified. We first treat the case of double
$\eta$-quotients, which is slightly more involved than $k \geq 3$.

\begin {theorem}
\label {th:ramifieddouble}
Under the assumptions of Theorem~\ref {th:singulardouble}, let $p_1 \neq
p_2$ be ramified in~$K$. Denoting by $\pf_i$ the ideal of the maximal
order~$\Oc_\Delta$ of~$K$ above $p_i$ and by $\sigma (\pf_i)$ the associated
Frobenius automorphism of $\Omegac / K$, we have
\begin {enumerate}
\item
$\wpq^e (\tau)^{\sigma (\pf_1)} = \legendre {p_1}{p_2'}^e \frac {1}{\wpq^e
(\tau)}$,
\item
$\wpq^e (\tau)^{\sigma (\pf_2)} = \legendre {p_2}{p_1'}^e \frac {1}{\wpq^e
(\tau)}$, and
\item
$\wpq^e (\tau)^{\sigma (\pf_1 \pf_2)} =
\left\{ \begin {array}{ll}
(-1)^{\frac {e (p_1 - 1)(p_2 - 1)}{4}} \, \wpq^e (\tau)
& \text { if } 2 \nmid p_1 p_2 \\
(-1)^{\frac {e (p_2^2 - 1)}{8}} \, \wpq^e (\tau)
& \text { if } p_1 = 2
\end {array} \right.$.
\end {enumerate}
In particular, if $|D| \not\in \{ p_1 p_2, 4 p_1 p_2 \}$ and one of the
following conditions holds:
\begin {itemize}
\item
$e$ is even;
\item
$p_1$ and $p_2$ are odd and one of them is congruent to~$1$ modulo~$4$;
\item
$p_1 = 2$ and $p_2 \equiv \pm 1 \pmod 8$;
\end {itemize}
then $\wpq^e (\tau)$ lies in the subfield of index~$2$ of $\Omegac / K$ with
Galois group $\Cl (\OD) / \langle \pf_1 \pf_2 \rangle$,
and $H_D^{\wpq^e}$ is the square of a polynomial in $K [X]$ (resp. $\Q [X]$
if $H_D^{\wpq^e} \in \Q [X]$).
\end {theorem}

\begin {proof}
We rely on Shimura's reciprocity law in the formulation of
\cite[Th.~5.1.2]{Schertz10}. Since $p_1$ and $p_2$ divide the level of
$\wpq^e$, we cannot use it directly; instead, we apply it twice for $\pf_1$
on the singular value $\wq^e (\tau)$. Recall from \S1 that $\wq^e$ is
modular of level $n = \frac {s}{\gcd (e, s)} \, p_2$ with $s = \frac
{24}{\gcd (24, p_2 - 1)}$; the hypotheses of Theorem~\ref
{th:singulardouble} imply that $e$ is sufficiently large so that $p_1$ is
coprime to the level of $\w1^e$. Since $N = p_1 p_2$ divides~$C$ and $p_1$
does not divide~$c$ by assumption, $\tau$ is the basis quotient of the ideal
$\af = \begin {pmatrix} \frac {-B + \sqrt D}{2} \\ A \end {pmatrix}_{\Z}$,
and $\frac {\tau}{p_1}$ is the basis quotient of $\af \bar \pf_1 = \af \pf_1
= \begin {pmatrix} \frac {-B + \sqrt D}{2} \\ p_1 A \end {pmatrix}_{\Z}$;
the matrix $P_1 = \begin {pmatrix} 1 & 0 \\ 0 & p_1 \end {pmatrix}$ of
determinant~$p_1$ sends the former to the latter basis. By Shimura
reciprocity, we have
\begin {eqnarray*}
\wq^e (\tau)^{\sigma (\pf_1)}
& = & \left( \wq^e \circ \left( p_1 P_1^{-1} \right) \right) (P_1 \tau)
= \left( \wp^e \circ S \circ P_1 \circ S \right) \left( \frac {\tau}{p_1}
\right) \\
& = & \left( \left( \sqrt {p_2} \, \frac {\eta (p_2 z)}{\eta (z)} \right)^e
\circ P_1 \circ S \right) \left( \frac {\tau}{p_1} \right).
\end {eqnarray*}
The action of $P_1$ on $\sqrt {p_2}$ is given by multiplication by $\xi \in
\{ \pm 1 \}$, and it is trivial on the rational $q$-expansion of $\frac
{\eta (p_2 z)}{\eta (z)}$. Thus,
\[
\wq^e (\tau)^{\sigma (\pf_1)} = \xi^e \wq^e \left( \frac {\tau}{p_1} \right)
= \xi^e \frac {\wq^e (\tau)}{\wpq^e (\tau)}.
\]
A second application of $\sigma (\pf_1)$ and using $\pf_1^2 = (p_1)$ and
$\xi^2 = 1$ yields
\begin {equation}
\label {eq:prooffrobenius}
\wq^e (\tau)^{\sigma (p_1)} = \frac {\wq^e (\tau)}{\wpq^e (\tau) \wpq^e
(\tau)^{\sigma (\pf_1)}}.
\end {equation}
The action of $\sigma (p_1)$ is again computed by Shimura reciprocity as
\[
\wq^e (\tau)^{\sigma (p_1)} = \left( \wq^e \circ \begin {pmatrix} p_1 & 0 \\
0 & p_1 \end {pmatrix} \right) (\tau).
\]
Notice that from the hypotheses of Theorem~\ref {th:singulardouble}, we have
$\gcd (p_1, n) = 1$. Bézout's relation between $p_1$ and $n^2$ yields a
matrix $M \in \Gamma$ with $M \equiv \begin {pmatrix} p_1 & 0 \\ 0 &
p_1^{-1} \end {pmatrix} \pmod n$, so that $\begin {pmatrix} p_1 & 0 \\ 0 &
p_1 \end {pmatrix} \equiv \begin {pmatrix} 1 & 0 \\ 0 & p_1^2 \end {pmatrix}
M \pmod n$. Then
\[
\wq^e \circ \begin {pmatrix} p_1 & 0 \\ 0 & p_1 \end {pmatrix} = \wq^e \circ
\begin {pmatrix} 1 & 0 \\ 0 & p_1^2 \end {pmatrix} \circ M = \wq^e \circ M =
\legendre {p_1}{p_2'}^e \wq^e
\text { by \eqref {eq:trafosimpleeta}}.
\]
Plugging this into \eqref {eq:prooffrobenius} finishes the proof of the
first formula of the theorem. The second formula is shown in the same way by
exchanging the roles of~$p_1$ and~$p_2$, and the third one follows from
quadratic reciprocity. Under the additional restrictions on $e$, $p_1$ and
$p_2$, it immediately follows that $\wpq^e (\tau)$ is invariant under
$\sigma (\pf_1 \pf_2)$.

The Galois automorphism $\sigma (\pf_1 \pf_2)$ is non-trivial unless
$(\pf_1 \pf_2) \cap \OD$ is principal.
This can be checked using genus theory of non-maximal orders,
or in an elementary fashion as follows: Write $|D| = p_1 p_2 r$ with $r \geq
1$. Then $(\pf_1 \pf_2) \cap \OD = \frac {t + v \sqrt D}{2} \, \OD$
if and only if
\[
p_1 p_2 = \frac {t^2 + v^2 |D|}{4} = \frac {t^2 + v^2 p_1 p_2 r}{4}
= p_1 p_2 \frac {s^2 p_1 p_2 + v^2 r}{4}
\]
with $t = s p_1 p_2$. This happens exactly for $s = 0$, $v = 2$, $r = 1$,
$|D| = p_1 p_2$; or $s = 0$, $v = 1$, $r = 4$, $|D| = 4 p_1 p_2$. So
excluding these cases, $\wpq^e$ lies in the subfield
of $\Omegac$ of index~$2$ and of Galois group $\Cl (\OD) / \langle \pf_1
\pf_2 \rangle$ over~$K$.
\end {proof}

For $k \geq 3$, the values of
$\w {p_1, \ldots, p_{i-1}, p_{i+1}, \ldots, p_k}^e (\tau)^{\sigma (\pf_i)}$
and
$\wpk (\tau)^{\sigma (\pf_i)}$ are computed in essentially the same way as
in the proof of Theorem~\ref {th:ramifieddouble}, but with $\xi$ and the
Legendre symbols dropped. This simplifies the argument and leads to fewer
restrictions on~$e$ and the~$p_i$:

\begin {theorem}
\label {th:ramifiedmulti}
Under the assumptions of Theorem~\ref {th:singularmulti}, in particular $k
\geq 3$, let $p_1, \ldots, p_k'$ for $2 \leq k' \leq k$ be ramified in $K$.
Denoting by $\pf_i$ the ideal of the maximal
order~$\Oc_\Delta$ of~$K$ above $p_i$ and by $\sigma (\pf_i)$ the associated
Frobenius automorphism of $\Omegac / K$, we have
\[
\wpk^e (\tau)^{\sigma (\pf_i)} = \frac {1}{\wpk^e (\tau)} \text { for } 1
\leq i \leq k'.
\]
If $k'$ is odd, or $|D| \not\in \{ p_1 \cdots p_{k'}, 4 p_1 \cdots p_{k'}
\}$, then $\wpk^e (\tau)$ lies in the subfield of index $2^{k'-1}$ of
$\Omegac$ which has Galois group $\Cl (\OD) / \langle \{ \pf_1 \pf_j : 2
\leq j \leq k' \} \rangle$ over $K$.
In particular, $H_D^{\wpk^e}$ is the $2^{k' - 1}$-th root of a polynomial in
$K [X]$ (resp. $\Q [X]$ if $H_D^{\wpk^e} \in \Q [X]$).
\end {theorem}

\begin {proof}
The action of $\sigma (\pf_i)$ is computed as above and implies that the
singular value is invariant under all the $\sigma (\pf_1 \pf_j)$ for $2 \leq
j \leq k'$. The classes of these ideals generate a subgroup of the
$2$-torsion subgroup of $\Cl (\OD)$. By the same argument as in the proof of
Theorem~\ref {th:ramifieddouble}, a product $\nf$ of several $\pf_1 \pf_j$
is principal if and only if $\nf = \pf_1 \cdots \pf_{k'}$, which can only
happen for even $k'$, and $|D| \in \{ \pf_1 \cdots \pf_{k'}, 4 \pf_1 \cdots
\pf_{k'} \}$.
\end {proof}

\begin {corollary}
\label {cor:unit}
Under the assumptions of Corollary~\ref {cor:realeven}, let $k' \geq 1$ of
the primes be ramified in $\OD$. Then $\wpk^e (\tau_1)$ is a unit.
For $k = 2$, if $p_i$ is ramified (and $p_{3-i}$ potentially not), then the
constant coefficient of $H_D^{\wpq}$ is $\legendre {p_i}{p_{3-i}'}^{\frac {e
h (\OD)}{2}}$.
Under the assumptions of Theorem~\ref {th:ramifieddouble}, the constant
coefficient of $\sqrt {H_D^{\wpq}}$ is $\legendre {p_1}{p_2'}^{\frac {e h
(\OD)}{4}}$.
For $k \geq 3$, the constant coefficient of $H_D^{\wpk^e}$ is $+1$.
Under the assumptions of Theorem~\ref {th:ramifiedmulti} and $|D| \not\in \{
p_1 \cdots p_{k'}, 4 p_1 \cdots p_{k'} \}$, the constant coefficient of
$\sqrt [2^{k'-1}]{H_D^{\wpk^e}}$ is $+1$.
\end {corollary}

\begin {proof}
Since $\af_i$ and $\af_i \pf_1$ are always inequivalent in $\Cl (\OD)$, each
conjugate $\wpk^e (\tau_i)$ occurs with its inverse (possibly up to sign if
$k=2$) by Theorems~\ref {th:ramifieddouble} and~\ref {th:ramifiedmulti}.
This shows the unit property and the results on the constant coefficients of
$H_D^{\wpq}$ and $H_D^{\wpk^e}$. The constant coefficient of the $2$-power
root has the desired property as long as $\pf_1$ does not lie in the
subgroup of $\Cl (\OD)$ generated by $\pf_1 \pf_2$ resp. the $\pf_1 \pf_2,
\ldots, \pf_1 \pf_{k'}$, which holds under the additional assumptions.
\end {proof}

\subsection {$n$-systems and reality of class polynomials, again}
\label {ssec:reality}

In the setting of Theorem~\ref {th:ramifieddouble}, the class polynomial is
the square of a polynomial with real coefficients, and we would like to
identify the elements of an $n$-system leading to identical values, thus
effectively cutting down the number of function evaluations to about $\frac
{h (\OD)}{4}$. We know by Theorem~\ref {th:ramifieddouble} that $\wpq^e
(\tau_{i_1}) = \wpq^e (\tau_{i_1})^{\sigma (\pf_1 \pf_2)}$. On the other
hand, the proof that an $n$-system yields the algebraic conjugates of
$\wpq^e (\tau_{i_1})$ relies on the equation $\wpq^e (\tau_{i_1})^{\sigma
(\pf_1 \pf_2)} = \wpq^e (\tau_{i_2})$, where $\tau_{i_2}$ is the basis
quotient of $\af_{i_2} \sim \af_{i_1} (\pf_1 \pf_2)^{-1} \sim \af_{i_1}
\nf$. So $\af_{i_1}$ and $\af_{i_2}$ determine a double root of the class
polynomial.

Using also Theorem~\ref {th:realnsystem}, we see that if $\af_{i_1}$,
$\af_{i_2} \sim \af_{i_1} \nf$, $\af_{i_3} \sim \af_{i_1}^{-1}$ and
$\af_{i_4} \sim \af_{i_1}^{-1} \nf$ correspond to four different elements of
the $n$-system, that is, if they are pairwise inequivalent, then they yield
twice the same pair of complex conjugate values.
If $\af_{i_1} \sim \af_{i_3}$, that is, $\af_{i_1}^2 \sim \OD$, then
$\af_{i_1}$ and $\af_{i_2} \sim \af_{i_1}^{-1} \nf$ yield the same real
value.
If $\af_{i_1} \sim \af_{i_4}$, that is, $\af_{i_1}^2 \sim \nf$, then again
$\af_{i_1}$ yields a real value by Theorem~\ref {th:realnsystem}, and
$\af_{i_2}$ yields the same real value.

This argumentation carries over immediately to Theorem~\ref
{th:ramifiedmulti}: For a given $\af_{i_1}$, the $2^{k' - 1}$ ideals of the
$n$-system that are equivalent to one of $\af_{i_1} (\pf_1 \pf_2)^{e_2}
\cdots (\pf_1 \pf_{k'})^{e_{k'}}$ with $e_2, \ldots, e_{k'} \in \{ 0, 1 \}$
lead to the same value of $\wpk^e$, so that the function needs to be
evaluated essentially only $\frac {h (\OD)}{2^{k' - 1}}$ times. If $k$ is
even, complex conjugate values may again by identified using Theorem~\ref
{th:realnsystem}, and the number of function evaluations drops to about
$\frac {h (\OD)}{2^{k'}}$.

Yet another factor of~$2$ may be saved by exploiting the action of $\sigma
(\pf_1)$. According to Theorems~\ref {th:ramifieddouble} and~\ref
{th:ramifiedmulti}, if $\af_{i_2} \sim \af_{i_1} \pf_1$, then $\wpk^e
(\tau_2) = \xi \frac {1}{\wpk^e (\tau_1)}$ with $\xi = 1$ for $k \geq 3$ and
$\xi = \legendre {p_1}{p_2'}^e$ for $k = 2$. This cuts down the number of
function evaluations to about $\frac {h (\OD)}{2^{k'}}$, or even $\frac {h
(\OD)}{2^{k' + 1}}$ when $k$ is even and thus the class polynomial is real.
In this optimal case, the $2^{k'+1}$ ideals of the $n$-system equivalent to
$\af_{i_1}^{e_0} \pf_1^{e_1} \cdots \pf_{k'}^{e_{k'}}$ with $e_0 \in \{ \pm
1 \}$ and $e_1, \ldots, e_{k'} \in \{ 0, 1 \}$ lead to $2^{k'-1}$ times the
same quadruple of two complex conjugate values and their (possibly negative)
inverses, assuming that the ideals are pairwise inequivalent, which happens
generically.

Even if only one of the primes is ramified, its explicit action gives useful
information.

\begin {corollary}
\label {cor:realodd}
Under the assumptions of Corollary~\ref {cor:realeven}, let $k$ be odd and
$p_1$ ramified in $\OD$. Then the class polynomial is real, and its complex
conjugate roots may be identified by Theorem~\ref {th:realnsystem} with $\cf
= \pf_1$.
\end {corollary}

\begin {proof}
Let $\tau_0 = A_1 \tau_1$ as in Theorem~\ref {th:realnsystem}. By \eqref
{eq:fal} and the action of $\sigma (\pf_1)$ according to Theorem~\ref
{th:ramifiedmulti} we have
\[
\wpk^e\alN (\tau_0) = \frac {1}{\wpk^e (\tau_0)}
= \wpk^e (\tau_0)^{\sigma (\pf_1)}.
\]
\end {proof}

\subsection {Examples}
\label {ssec:examplesramified}

Consider first $D = -455 = -5 \cdot 7 \cdot 13$. To simplify the
presentation, we use in the following the notation of quadratic forms,
identifying the ideal $\begin {pmatrix} \frac {-B + \sqrt D}{2} \\ A \end
{pmatrix}_{\Z}$ with the quadratic form $[A, B, C]$, where $C = \frac {B^2 -
D}{4 A}$.  Let $\pf_1 = [5, 5, 24]$, $\pf_2 = [7, 7, 18]$ and $\zf = [2, 1,
57]$. Then $\Cl (\OD) = \langle \pf_1, \zf \rangle \simeq \Z / 2 \Z \times
\Z / 10 \Z$ with $\zf$ of order~$10$ and $\pf_2 \sim \zf^5$.

The function $\w {5,7}$ satisfies the conditions of Theorem~\ref
{th:ramifieddouble} and is of level $n = 35$. A first element of a
$35$-system with $C$ divisible by~$35$ is given by $\af_1 = [1, 35, 420]
\sim \OD$. Then $\tau_1 = -17.5 + 10.66\ldots i$, and $z_1 = \w {5, 7}
(\tau_1) = -6.02\ldots$ is real. Since $\legendre {5}{7} = -1$, we have
another real conjugate $\frac {-1}{z_1}$ corresponding to $\pf_1$, and the
same two values reoccur for $\pf_1 \pf_2 \sim \pf_1 \zf^5$ and $\pf_1 (\pf_1
\pf_2) \sim \zf^5$.

The not yet covered ideal $\zf$ is equivalent to $\af_2 = [2, 105, 1435]$ in
the same $35$-system, yielding $\tau_2 = -26.25 + 5.33\ldots i$ and $z_2 =
\w {5,7} (\tau_2) = 0.65\ldots - 2.05\ldots i$. We then obtain twice the
values $z_2$, $\bar z_2$, $\frac {-1}{z_2}$ and $\frac {-1}{\bar z_2}$,
corresponding to the ideals $\zf$, $\pf_1 \zf^4$, $\pf_1 \zf$, $\zf^4$, and
the second time to $\pf_1 \zf^6$, $\zf^9$, $\zf^6$, $\pf_1 \zf^9$,
respectively.

Finally, $\pf_1 \zf^2$ has not occurred yet; it is equivalent to $\af_3 =
[3, 175, 2590]$ in the $35$-system with $\tau_3 = -29.16\ldots + 3.55\ldots
i$ and yields twice the conjugates $z_3 = \w {5,7} (\tau_3) = 1.50\ldots -
0.53\ldots i$, $\bar z_3$, $\frac {-1}{z_3}$ and $\frac {-1}{\bar z_3}$.

Let $x_i = \Re (z_i)$ and $n_i = z_i \bar z_i$, and define $g_i$ to be the
polynomial with coefficients in $\R$ whose roots are the conjugates related
to $z_i$. Then $g_1 = X^2 - \left( z_1 - \frac {1}{z_1} \right) X - 1$ and
$g_i = X^4 + 2 x_i \left( \frac {1}{n_i} - 1 \right) (X^3 - X) + \left( n_i
+ \frac {1 - 4 x_i^2}{n_i} \right) + 1$ for $i \in \{ 2, 3 \}$. Multiplying
$g_1$, $g_2$ and $g_3$ and rounding the resulting coefficients to integers,
we obtain
\[
\sqrt {H_{-455}^{\w {5,7}}} = X^{10} + 3 X^9 - 12 X^8 + 32 X^7 - 38 X^6 - 17
X^5 + 38 X^4 + 32 X^3 + 12 X^2 + 3 X - 1.
\]
Notice the constant coefficient $-1$ as predicted by Corollary~\ref
{cor:unit}.

As an example with more ramified primes, let us consider $D = -3795 = -3
\cdot 5 \cdot 11 \cdot 23$, which is divisible by~$3$, and all other prime
factors of which are congruent to~$-1$ modulo~$3$. Its class group is
isomorphic to $\Z / 2 \Z \times \Z / 2 \Z \times \Z / 4 \Z$. If we wish to
use a multiple $\eta$-quotient whose level is composed of only ramified
primes, then Theorems~\ref {th:singulardouble} and~\ref {th:singularmulti}
imply that we need to add an exponent of~$3$. By including three of the four
prime factors of~$D$, we may then gain a factor of~$4$ in the degree. For
instance, one computes
\[
\sqrt [4]{H_{-3795}^{\w {3, 5, 11}^3}} = X^4 - 200596 X^3 - 511194 X^2 -
200596 X + 1.
\]
Alternatively, we may add a split prime congruent to~$+1$ modulo~$3$ to the
level of the $\eta$-quotient, which enables us to drop the exponent~$3$.
This will generally result in a smaller polynomial:
\[
\sqrt [4]{H_{-3795}^{\w {3, 5, 11, 19}}} = X^4 - 46 X^3 + 2115 X^2 - 46 X +
1.
\]
It should be noted that the singular values of the multiple $\eta$-quotients
do not necessarily \textit {generate} the subfields indicated in
Theorems~\ref {th:ramifieddouble} and~\ref {th:ramifiedmulti}. Using the
split prime~$13$ instead of~$19$ above yields
\[
\sqrt [4]{H_{-3795}^{\w {3, 5, 11, 13}}} = X^4 + 92 X^3 + 2118 X^2 + 92 X +
1 = (X^2 + 46 X + 1)^2;
\]
apparently, the prime ideals of $\OD$ above~$13$ and~$19$ act differently on
the respective $\eta$-quotient, although they are all of order~$4$ in the
class group.

\section {Other functions on $X_0^+ (N)$ in the totally ramified case}

Shimura reciprocity relates the Galois action on singular values to
actions on modular functions and their arguments.
In \S\ref {sec:ramified} we used the precise shape of the multiple
$\eta$-quotients to show invariance results when two or more of the
prime factors are ramified.
We may consider in more generality functions on $X_0^+ (N)$, that is,
functions invariant under $\Gamma^0 (N)$ and the Fricke--Atkin--Lehner
involution $\alN$.
Assuming less knowledge about the function, we will make the stronger
hypothesis on the imaginary-quadratic field that \textit {all} prime
divisors of~$N$ are ramified and consider the image under the Artin map
of the product of all ramified prime ideals.
Under similar technical conditions as before, one can then show that the
singular values lie in a subfield of index~$2$ of the ring class field.

\begin {theorem}
\label {th:atkin}
Let $D$ be a quadratic discriminant of conductor $c$ and let $N$ be
square-free, prime to $c$ and such that $N \mid D$ and
$|D| \not\in \{ N, 4 N \}.$
Then there is a primitive
quadratic form $[A_1, B_1, C_1]$ of discriminant~$D$ such that $C_1 = N$ and
$N \mid B_1$. Let $\nf = \begin {pmatrix} \frac {-B_1 + \sqrt D}{2} \\ A_1
\end {pmatrix}$ be the ideal of basis quotient $\tau_1$ associated to this
quadratic form. Then $\nf$ is of order~$2$ in $\Cl (\OD)$.

Let $f$ be a modular function for $\Gamma^0 (N)$ such that $f$ and $f \circ
S$ have a rational $q$-expansion, and such that $f\alN = f$. Then the
singular value $f (\tau_1)$ lies in the subfield of index~$2$ of $\Omegac$
with Galois group $\Cl (\OD) / \langle \nf \rangle$. The class polynomial
$\sqrt {H_D^f}$ is real and can be computed from an $N$-system in which
$\af_i$ and $\af_i \nf$ yield the same root, and $\af_i^{-1}$ and
$\af_i^{-1} \nf$ yield the complex conjugate root.
\end {theorem}

\begin {proof}
Since all primes dividing $N$ are not inert, there is a quadratic form with
$C_1 = N$. As $N$ divides $D$ and is square-free, it follows that $N \mid
B_1$ and that all primes dividing~$N$ are ramified. The ideal $\nf$ is
equivalent to the ideal of norm~$N$; it is non-principal by the same
argument as in the proof of Theorem~\ref {th:ramifieddouble}.

\cite[Th.~4]{Schertz02}, given in the beginning of \S\ref {ssec:classinv},
implies that $f (\tau_1) \in \Omegac$. Let $[A_i, B_i, C_i]$ be an
$N$-system, associated to the ideals $\af_i$ of basis quotients $\tau_i$.
The class polynomial $H_D^f$ is real by Theorem~\ref {th:realnsystem};
precisely, the ideals $\af_i$ and $\af_i^{-1} \nf$ yield complex conjugate
roots. Now,
\[
f (\tau_i)
= f\alN (\tau_i)
= f \left( \frac {-N}{\tau_i} \right)
= f \left( \frac {B_i + \sqrt D}{2 C_i / N} \right)
= f (\tau_j),
\]
where $\tau_j$ is a root of the quadratic form $[C_i / N, -B_i, A_i N]$ of
discriminant~$D$. From the $N$-system condition and $N \mid B_1$ we have $N
\mid B_i$ and $-B_i \equiv B_i \pmod {2 N}$. Moreover, $\gcd (C_i / N, N) =
1$ since $N \mid B_i$, $N \mid C_i$ and $N$ is coprime with~$c$. So the form
is primitive and can be assumed to occur as an element of the $N$-system,
associated to an ideal $\af_j \sim \af_i \nf$. Since $\nf$ is not principal,
$\af_i$ and $\af_j$ are not equivalent and correspond to different elements
of the $N$-system.
\end {proof}

Applied to multiple $\eta$-quotients, Theorem~\ref {th:atkin} is weaker
than Theorem~\ref {th:ramifieddouble} (which also considers powers of
$\wpq$ of lower exponent $e < s$, for which the function is of level
$\frac {s}{e} N$ instead of~$N$) and Theorem~\ref {th:ramifiedmulti}
(which treats the case that only some primes are ramified).
Its interest lies in its application to other functions invariant under the
Fricke--Atkin--Lehner involution. In particular, it was suggested in \cite
{Morain09} to consider such functions~$A_p$ of prime level~$p$ and of
(conjectured) minimal degree, as are used inside the Schoof--Elkies--Atkin
algorithm for point counting on elliptic curves \cite{Elkies98,Morain95}. It
was observed in \cite{Morain09} that $A_{71}$ yields an asymptotic gain in
the logarithmic height of the class polynomial, compared to the polynomial
for the $j$-invariant, by a factor of~$36$ as $|D| \to \infty$; families of
functions with a conjectured asymptotic gain of at least a factor of~$30$
are given in \cite{EnSu10}. Using ramified primes has a double advantage:
Not only is the degree of the class polynomial divided by~$2$, but also the
coefficients of the polynomial have about half as many 
digits.

\section {Class field theoretic remarks}

The main algorithmic use of class polynomials is the computation of elliptic
curves over some finite field~$\F$ with a given endomorphism ring~$\OD$. If
such curves exist, the class polynomial $H_D^f$ splits completely in $\F
[X]$; for each root $\bar f$, there is an elliptic curve over $\F$ with
complex multiplication by~$\OD$. To compute such a curve, one considers the
modular polynomial $\Psi_f (f, j)$, the minimal polynomial of $f$ over $\C
(j)$, which is in fact an element of $\Z [f, j]$ for all functions
considered in this article. So one may reduce the polynomial to an element
of $\F [f, j]$, specialise $f$ as $\bar f$ and compute all the roots $\bar j
\in \F$. For all these roots, of which there may be several, one writes down
an elliptic curve over $\F$ with $j$-invariant $\bar j$ and checks whether
its endomorphism ring is indeed $\OD$ as desired. For this application, it
is clearly computationally advantageous to compute only $\sqrt [2^{k'-1}]
{H_D^f}$.

If one wishes to obtain the full class field $\Omegac$, one may use the
subfield $L = K (x) = K [X] / \left( \sqrt [2^{k'-1}] {H_D^f (X)} \right)$
(assuming the polynomial is irreducible, which need not be the case, see the
example in \S\ref {ssec:examplesramified}), as a first step in a tower of
field extensions. To this purpose, one may factor the modular polynomial
$\Psi_f (x, Y)$ over~$L$. It necessarily has an irreducible factor of
degree~$2^{k'-1}$, which generates $\Omegac / L$. In general, this
polynomial factorisation step over a number field will be more costly than
switching to a different class invariant that directly generates $\Omegac$.
For a few functions of low level, however, the degree $d_j$ of $\Psi_f$
in~$j$ is exactly $2^{k-1}$, and no factorisation is needed for a
discriminant for which all primes dividing~$N$ are ramified. For $k = 2$, a
degree $d_j = 2$ is obtained if and only if $(p_1 - 1)(p_2 - 1) \mid 24$ by
\cite[Th.~9]{EnSc05}, that is, for $\{ p_1, p_2\} \in \big\{ \{2, 3\}, \{2,
 5\}, \{2, 7\}, \{2, 13\},  \{3, 5\}, \{3, 7\}, \{3, 13\}, \{5, 7\} \big\}$.
Conjecturally, for $k \geq 3$ we have $d_j = 2^{k-1}$ if and only if $(p_1 -
1) \cdots (p_k - 1) \mid 24$, that is, for $\{ p_1, \ldots, p_k \} \in
\big\{ \{ 2, 3, 5 \}, \{ 2, 3, 7 \}, \{ 2, 3, 13 \}, \{ 2, 5, 7 \} \big\}$.
Computing the polynomials via the algorithm of \cite{Enge09mod} shows that
the condition is at least sufficient, and that the smallest function with $k
= 4$, $\w{2,3,5,7}$, has $d_j = 16 > 2^3$. A top-down approach to obtain the
class field as a tower of field extensions, starting with a generating
element of $\Omegac$, is described in \cite{HaMo01,EnMo03}. The lucky case
$d_j = 2^{k-1}$ can be seen as a bottom-up approach, in which moreover the
second stage is realised by the universal modular polynomial independently
of the concrete class field under consideration.

It may also be possible to construct the missing part of $\Omegac$
classically using genus theory. Suppose that $\sqrt [2^{k'-1}] {H_D^f (X)}$
is irreducible. Let $H = \langle \{ \pf_1 \pf_j : 2 \leq j \leq k' \}
\rangle$, and assume that $\Cl (\OD)$ is the direct product of $\Cl (\OD) /
H$ and $H$, which happens if and only if $H$ contains no element that is a
square in $\Cl (\OD)$. Then $\Omegac$ is the compositum of~$L$ and the genus
field of~$K$.

For instance, consider the first example of \S\ref {ssec:examplesramified},
where $D = - 5 \cdot 7 \cdot 13$ and $L / K$ is generated by $\sqrt
{H_{-455}^{\w{5,7}}}$. A quick computation with Pari/GP \cite {parigp} shows
that $L \ni \sqrt {13}$, but $L \not\ni \sqrt 5$, so that $\Omegac = L
(\sqrt 5) = L (\sqrt {-7})$.

However, the singular values of the multiple $\eta$-quotients do not
necessarily generate the Hilbert class field of a fundamental discriminant
over the genus field. For $D = -3 \cdot 5 \cdot 11 \cdot 23$ and $L / K$
generated by $\sqrt [4]{H_{-3795}^{\w {3, 5, 11, 19}}}$, one has $L \ni
\sqrt {-3}, \sqrt {-11}$, and $L (\sqrt 5) = L (\sqrt {-23})$ still has
index~$2$ in~$\Omegac$.

\bibliographystyle {alpha}
\bibliography {zweig}

\end{document}